%% file: cattiauxguillinmalrieu.tex
\documentclass [a4paper,twoside,11pt]{article}

\usepackage[latin1]{inputenc}
\usepackage{amsmath,amsfonts,amssymb}
\usepackage{vmargin,graphicx,theorem}
\usepackage[english]{babel}
\usepackage{enumerate}

\setpapersize[portrait]{A4}
\setmarginsrb{2.2cm}{2cm}{2.2cm}{1cm}{0cm}{0.1cm}{0.5cm}{2cm}

\include{macros}

\begin{document}

\title{\sl Probabilistic Approach for Granular Media Equations in the Non Uniformly Convex Case}
\author{P. Cattiaux, A. Guillin and F. Malrieu }

\date{\today}
\maketitle\thispagestyle{empty}

\begin{abstract}
  We use here a particle system to prove a convergence result as well
  as a deviation inequality for solutions of granular media equation
  when the confinement potential and the interaction potential are no
  more uniformly convex. Proof is straightforward, simplifying deeply
  proofs of Carrillo-McCann-Villani \cite{CMV,CMV2} and completing results
  of Malrieu \cite{malrieu03} in the uniformly convex case. It relies
  on an uniform propagation of chaos property and a direct control in
  Wasserstein distance of solutions starting with different initial
  measures. The deviation inequality is obtained via a $T_1$
  transportation cost inequality replacing the logarithmic Sobolev
  inequality which is no more clearly dimension free.
\end{abstract}

\bigskip

{\it Mathematics Subject Classification 2000:} 65C35, 35K55, 65C05, 82C22, 26D10, 60E15.\\
{\it Keywords:} Granular media equation, transportation cost
inequality, Logarithmic Sobolev Inequalities - Concentration
inequalities.

\par\vspace{40pt}


\section{Introduction}

Our main goal will be to deal in a probabilistic way with the
following nonlinear equation
\begin{equation}
  \label{eq:edp}
\PD{u}{t}=\nabla\cdot[\nabla u+u\nabla V+u\nabla W *u]
\end{equation}
where $u(t,x)$ is a time dependent probability measure, $*$ denotes
the standard convolution operator and $V$ and $W$ are two convex (at
infinity) potentials. \par\vspace{5pt}

This equation arises (in dimension 1) in the modeling of granular
media as follows: consider many infinitesimal particles colliding
inelastically. With a correct renormalization between the frequency
and the inelasticity of the collisions, $u(t,x)$ turns out to be the
velocity of a representative particle (among an infinity). The
potential $V$ represents the friction and $W$ the inelastic collisions
between particles with different velocities. Note that the particular
case $V=0$ and $W(x)=|x|^3$ is of special interest and we refer to
Benedetto-Caglioti-Pulvirenti \cite{BCP} for the physical issues (also
see \cite{BCCP}).\par\vspace{5pt}

Once the problem of existence and uniqueness is tackled, one major problem in this equation is the
behavior at infinity: existence of a stationary measure and speed of convergence towards this
stationary measure or even distance between two solutions starting at different points. It has
been studied in parallel by Carrillo-McCann-Villani \cite{CMV,CMV2} and Malrieu
\cite{malrieu01,malrieu03}, under various assumptions on the potentials $V$ and $W$, using
analytical and probabilistic approaches respectively (also see \cite{brv98,brtv98} for one
dimensional particles). We will consider here the probabilistic approach and will recover and
generalize slightly results of \cite{CMV,CMV2} as well as give a quantified probabilistic
approximation of the stationary measure of the granular media equation. It is worthwhile noticing
that the analytic methods in \cite{CMV,CMV2} cover a much larger spectrum of non linear p.d.e's
(like the porous medium equation), for which the probabilistic approach remains to be written.

As physical interest (and in fact where the main mathematical
difficulty resides) dictates the friction term to vanish, we will
consider the following two sets of general assumptions:
\begin{edefi}\label{assump}
\begin{itemize} \item  We say that $W$ satisfies the set of assumptions (A) if
\begin{enumerate}
\item[A1.] the friction term $V=0$; \item[A2.] $W$ is symmetric, i.e.
  $W(-x)=W(x)$; \item[A3.]  $\nabla W$ is locally Lipschitz with
  polynomial growth, i.e. for some $m$
  $$
  \forall x,y\in\dR^d, |\nabla W(x)-\nabla W(y)|\le C (|x-y|\wedge
  1) (1+|x|^m+|y|^m).
  $$
  We also assume that the second derivative of $W$ has a ($m$)
  polynomial growth.
\item[A4.] $W$ is the sum of a compactly supported $\mathcal{C}^2$
  function and a $\mathcal{C}^2$ uniformly convex function.

  This last assumption entails some uniform convexity at infinity
  property, namely that there exist \underline{positive} $C, \lambda$,
  such that
\begin{equation}\label{cond:conv}
\forall x,y,\qquad (x-y).(\nabla W(x)-\nabla W(y))\ge \lambda \|x-y\|^2-C \, .
\end{equation}
\end{enumerate}

\item We say that $V$ and $W$ satisfy the set of assumptions (A') if
  $V$ is uniformly convex at infinity (i.e. satisfies (A4)), $W$ is
  convex at infinity (but not necessarily uniformly) and symmetric,
  and both satisfy (A3).
\end{itemize}
\end{edefi}

As we will see later, these assumptions are not sufficient to get good
properties for large time of the granular media equation and will be
replaced by some ``strict convexity except for a finite number of
points'' property.

\par\vspace{5pt}

Before further discussing the result, let us first present the
probabilistic approach of this problem. The probabilistic
interpretation is to consider a Markov process $(\bar X_t)_{t\ge 0}$,
which law at time $t$ is $u$. It is the solution of the nonlinear S.D.E.
\begin{equation}
  \label{eq:nonlin}
\left\{\begin{array}{l}
d\bar X_t=\sqrt{2}dB_t-\nabla V(\bar X_t)dt-\nabla W*u_t(\bar X_t)dt\\
{\cal L}(\bar X_t)=u_t \, dx.\end{array}\right.
\end{equation}

We wish here first to give sufficient conditions ensuring both
existence and ergodicity of the solution of the nonlinear S.D.E. and in
a second time to provide a way to simulate this law at each time $t$
with some Gaussian confidence intervals independent of time $t$. These
two goals will be carried out through the extensive use of some
(linear) particle approximations, i.e. $(X^N_t)$ solution of
\begin{equation}
  \label{eq:particle}
\left\{\begin{array}{ll}
dX^{i,N}_t=\sqrt{2}dB^i_t-\nabla V(X^{i,N}_t)dt-{1\over N}\sum_{j=1}^N\nabla W(X^{i,N}_t-X^{j,N}_t)dt&i=1,...,N\\
X^{i,N}_0=X^i_0&i=1,...,N.\end{array}\right.
\end{equation}

If we suppose that $V=0$, we may assume that the center of mass is
fixed, and without loss of generality, set to $0$, i.e. $\dE \bar
X_0=0$, indeed it is easy to remark that $\bar X_t-\dE\bar X_0$
satisfies equation (\ref{eq:nonlin}). This assumption obliges to
introduce the projected system onto the set $\sum_i x_i=0$ (see
Section \ref{2}), and introduces some intricacies.  \smallskip

Under the strong assumptions of (global) uniform convexity of one of the potentials (the other one
being globally convex but non necessarily uniformly), Malrieu \cite{malrieu01,malrieu03}
successively fulfilled these two goals using extensively the so called Bakry-\'Emery criterion
ensuring that a Logarithmic Sobolev inequality holds independently of $N$. He also provided in the
same time asymptotic behavior in time and concentration inequality, enabling him to recover part
of the results of Carrillo-McCann-Villani \cite{CMV} (see also Bolley-Guillin-Villani \cite{BGV}
for a strengthened deviation inequality for almost quadratic potentials).

The assumption of uniform convexity is however far too strong for applications (preventing for
example to consider the case of cubic interaction potential). It can be removed (see \cite{CMV})
for a one point degeneracy where the authors obtain various rates of convergence. Another approach
is proposed by Carrillo-McCann-Villani in the subsequent \cite{CMV2} paper introducing tools of
contractions in $L^2$-Wasserstein distance for length space enabling them to ensure a convergence
in Wasserstein distance of the solution of the granular media equation.
\smallskip

We will see here that, with much simpler tools, we can recover their
result and may provide also Gaussian confidence bounds for the
approximations of the granular media equations. It is worthwhile
noticing that, though the invariant measures of the particle system
still satisfy some log-Sobolev inequality, it is difficult to obtain
dimension-free estimates on the log-Sobolev constants (as for nearest
neighbors interactions models in Statistical Mechanics). This prevents
us to use this kind of approach.  \smallskip

Our main tool in order to get the Gaussian concentration inequality we
need, will thus be a transportation cost-information inequality: let
say that $\mu\in T_1(C)$ if for every probability measure $\nu$
\begin{equation}
  \label{eq:transp}
W_1(\nu,\mu)\le\sqrt{C \, \text{Ent}(\nu|\mu)}
\end{equation}
where $W_1$ is the usual Wasserstein distance and $Ent$ the Kullback
information or relative entropy. In the final section we prove that,
under one of the hypothesis (A) or (A'), the law of the particle
system at time $t$ satisfies a $T_1(C N)$ inequality, for some $C$
that does not depend on $N$ nor on $t$.  \smallskip

To complete the proof of the concentration inequality in the final
section we need some ``uniform (in time) propagation of chaos''. It
seems difficult to obtain such a result only assuming (A) or (A').
Indeed if the potentials are non attractive in some (bounded) region,
the situation becomes the classical mean-field one where propagation
of chaos is controlled on finite time intervals only. That is why we
have to reinforce our assumptions.

Our condition below is inspired by the work of Carrillo-McCann-Villani
\cite{CMV2}:

we say that condition $\mathbf {C(A,\alpha)}$ holds if there exist
$A,\alpha>0$ such that for any $0<\epsilon<1$,
\begin{equation}\label{supercond}
\forall x,y\in\dR^d,\qquad (x-y). (\nabla W(x)-\nabla W(y))\ge A\epsilon^\alpha(
|x-y|^2-\epsilon^2).
\end{equation}

Remark that this condition implies convexity (if $|x-y|\leq 1/2$
choose $\epsilon=|x-y|$, otherwise choose $\epsilon=1/4$) but is
weaker than uniform convexity (which is true with $\alpha=0$). Typical
examples are polynomial potentials of the form $W(x)=|x|^{2+\alpha}$
which satisfy $\mathbf{ C(A,\alpha)}$ for some positive $A$. In fact
this condition is related to the rate at which strict convexity is
lost at points.

Condition $\mathbf {C(A,\alpha)}$ allows us to prove the required
``uniform (in time) propagation of chaos'' and then to obtain the rate
of convergence of $u_t$ to the limit $u_\infty$ in Wasserstein
distance $W_2$. These results are obtained in sections 3 and 4.
\smallskip

To complete the description of the paper, the next section contains
some useful estimates, and a complete proof of existence and
uniqueness for all the equations we have considered (\eqref{eq:edp}
and \eqref{eq:nonlin} in particular). The proof is quite natural
(using the particle system). We did not find such a proof in the
literature except in dimension 1, and as a matter of fact the set of
hypotheses we need, slightly differs from the ones used in the
analytical literature.  \bigskip


\section{The particle system and solutions of the nonlinear SDE}\label{2}

In this section we study the particle system and the nonlinear SDE.
The first discussion is about case (A) when $V$ vanishes.

Consider once again the particle system $(X^N_t)$ solution of
\begin{equation}\label{eq:part}
\left\{\begin{array}{ll}
dX^{i,N}_t=\sqrt{2}dB^i_t - \nabla V(X^{i,N}_t) dt -{1\over N}\sum_{j=1}^N\nabla W(X^{i,N}_t-X^{j,N}_t)dt&i=1,...,N\\
X^{i,N}_0=X^i_0&i=1,...,N.\end{array}\right.
\end{equation}
Here $X^{i,N}_0$ are i.i.d. variables, with common law $\mu_0$.

As noted by Malrieu, when $V=0$, the direction $\{(v,\cdot,,v)\}$ is
quite singular for the particle system, and we then introduce its
orthogonal hyperplane ${\cal M}=\left\{x\in(\dR^d)^N;\sum_{i=1}^N
  x_i=0\right\}$, and consider the projected particle system
$$
Y^{i,N}_t=X^{i,N}_t-{1\over N}\sum_{j=1}^N X^{j,N}_t\qquad \forall
i=1,...,N.
$$
The process ${(Y^{N}_t)}_{t\geq 0}$ thus verifies the system of SDE's
\begin{equation}\label{eq:projpart}
\left\{\begin{array}{ll}
dY^{i,N}_t=\sqrt{2}dB^i_t-{\sqrt{2}\over N}\sum_{j=1}^NdB^j_t-{1\over N}\sum_{j=1}^N\nabla W(Y^{i,N}_t-Y^{j,N}_t)dt
&i=1,...,N\\
Y^{i,N}_0=X^i_0-{1\over N}\sum_{j=1}^N X^{j,N}_0&i=1,...,N\end{array}\right.
\end{equation}
hence is a diffusion on $\{\cal M\}$. Note that $\dE(Y^{i,N}_t)=0$.\par\vspace{10pt}

We assume that $W$ is satisfying the assumptions (A) in Definition
\ref{assump}. Let us first remark that this new particle system
satisfies globally the ``convexity at infinity'' property. We have to
verify, denoting $b(x)=(b^1(x),...,b^N(x))$, $b^i$ with values in
$\dR^d$ and
$$
b^i(x)=-\frac{1}{N}\sum_{j=1}^N\nabla W(x_i-x_j),
$$
that condition (\ref{cond:conv}) holds true.  Note that, with
$\sigma$ the (constant) matrix diffusion coefficient of
$(Y^{i,N}_t)_{1_le i\le N}$, we have $\|\sigma\|_{HS}\le c$ for some
positive $c$, and, as $W$ is symmetric and $x,y\in{\cal M}$
\begin{eqnarray*}
(x-y)\cdot(b(x)-b(y))&=&-\frac{1}{2N}\sum_{i,j=1}^N(x_i-x_j-(y_i-y_j))\cdot(\nabla W(x_i-x_j)-\nabla W(y_i-y_j))\\
&\le& \frac{1}{2N}\sum_{i,j=1}^N\left(-\lambda|x_i-x_j-y_i+y_j|^2+C\right)\\
&=&-\lambda \|x-y\|^2+\frac{CN}{2},
\end{eqnarray*}
where the latter equality follows as $x,y\in{\cal M}$.
\smallskip

This remark allows to prove existence, uniqueness and non-explosion
for the solution of the SDE \eqref{eq:projpart}. Indeed recall that a
sufficient condition for all this to hold is the following:

there exists some $\psi$ such that $\psi(x) \rightarrow +\infty$ as
$|x| \rightarrow +\infty$ and $\Delta \, \psi \, + \, b.\nabla \psi$
is bounded from above.

$\psi$ is a kind of Lyapounov function. Here we may choose
$\psi(x)=|x|^2$ according to the previous remark applied with $y=0$.
In this case non-explosion holds as soon as $\mu_0(|x|^2)<+\infty$.

The same argument can be used with (A') if $W$ is everywhere convex.

For the initial system \eqref{eq:part} the only thing to show is non
explosion. It will follow from Proposition \ref{prop:xmoment} below.
\medskip

\subsection{Moment controls for the particle system}

In order to prove tightness of the empirical measure of the particle
system, and for later use, we first prove some controls on moments.

\begin{eprop}\label{prop:xmoment}
  If $W$ and $V$ satisfy (A'), and $\mu_0(|x|^2)<+\infty$ the solution
  to \eqref{eq:part} is non explosive. Furthermore there exists some
  $K>0$ such that for all $i$,
  $$
  \sup_{t\geq 0} \dE|X_t^{i,N}|^2\le \mu_0(|x|^2) + K \, .
  $$
  If $W$ satisfies (A) a similar result holds with $Y_t^{i,N}$ in
  place of $X_t^{i,N}$.
\end{eprop}
\begin{eproof}
  Using It\^o's formula up to the infimum between $t$ and the exit
  time of a large ball (let $T$ be this infimum) we have
\begin{eqnarray*}
\dE\sum_{i=1}^N |X_T^{i,N}|^2&=&\dE\sum_{i=1}^N |X_0^{i,N}|^2+2dN \dE(T) - 2 \dE \sum_{i=1}^N
\int_0^T
X_s^{i,N} . \nabla V(X_s^{i,N}) ds \\
 & & -{2\over N}\dE\sum_{i,j=1}^N \int_0^T X_s^{i,N}.\nabla W(X_s^{i,N}-X_s^{j,N})ds\\
&=&\dE\sum_{i=1}^N |X_0^{i,N}|^2+2dN \dE(T)  - 2 \dE \sum_{i=1}^N \int_0^T X_s^{i,N} . \nabla V(X_s^{i,N}) ds\\
 & & -{1\over N}\dE\sum_{i,j=1}^N \int_0^T (X_s^{i,N}-X_s^{j,N}).\nabla W(X_s^{i,N}-X_s^{j,N})ds.
\end{eqnarray*}
We may use condition \eqref{cond:conv}, go to the limit with respect
to the radius of the ball (hence replace $T$ by $t$) and obtain the
finiteness of the quantity. This implies non explosion.

Furthermore, denoting $v(t)=\dE\sum_{i=1}^N |X_t^{i,N}|^2$,
differentiating (see the proof of the next Proposition) and using our
conditions (in particular $x.\nabla W(x) \geq \lambda_W |x|^2 - C$),
we get
\begin{eqnarray*}
v'(t)&\le& -2\lambda_V v(t)+2dN+3CN.
\end{eqnarray*}
Gronwall's lemma and exchangeability conclude the proof in case (A').

The proof is similar for $Y$ in case (A) using the convexity at
infinity of the drift, as we previously remarked.
\end{eproof}

\begin{eprop}\label{prop:momdiff}
  If $W$ satisfies (A) or $W$ and $V$ satisfy (A'), for all $k\in \dN$
  there exists $C(k)>0$ such that for all $1\leq i,j\leq N$,
  $$
  \sup_{t\geq 0} \dE\PAR{|X_t^{i,N}-X_t^{j,N}|^{2k}}
  \leq C(k)\left(1+\mu_0\otimes \mu_0\left(|x-y|^{2k}\right)\right).
  $$
\end{eprop}

\begin{eproof}
  We write the proof in the case $V=0$ (i.e (A)), the case (A') is
  similar.

  Recall that all particles are exchangeable. We may apply It\^o's
  formula up to the exit time of a large ball (for the whole system)
  and then go to the limit in order to get,
\begin{eqnarray*}
& & \dE\sum_{i,j=1}^N |X_t^{i,N}-X_t^{j,N}|^{2k} = \dE\sum_{i,j=1}^N |X_0^{i,N}-X_0^{j,N}|^{2k}+
2k(2k-1) \int_0^t \dE\sum_{i,j=1}^N |X_s^{i,N}-X_s^{j,N}|^{2k-2} ds  \\ & & - {2k\over N} \,
\dE\sum_{i,j,l=1}^N \int_0^t \left(\nabla W(X_s^{i,N}-X_s^{l,N})- \nabla
W(X_s^{i,N}-X_s^{l,N})\right) .(X_s^{j,N}-X_s^{l,N}) |X_s^{i,N}-X_s^{j,N}|^{2k-2}  ds.
\end{eqnarray*}
Denoting $A_k(t)=\dE\sum_{i,j=1}^N |X_t^{i,N}-X_t^{j,N}|^{2k}$, and
using \eqref{cond:conv} we obtain
$$
A_k(t)\leq N \mu_0\otimes \mu_0\left(|x-y|^{2k}\right) + 2k(2k-1+C)
\, \int_0^t A_{k-1}(s) ds - 2k \lambda \int_0^t A_k(s) ds \, .
$$
Applying Gronwall's lemma and an easy induction we thus have that
$A_k(t)$ is finite for all $t$.  Accordingly we may replace the pair
of times $(0,t)$ by $(t,t+\varepsilon)$ and prove that $t \mapsto
A_k(t)$ is differentiable.  Differentiating at time $t$ yields
$$A'_1(t) \leq -2 \lambda A_1(t) + 2N(d+C) \, .$$
Gronwall's lemma
yields the desired result for $k=1$. The proof follows, using this
bound and an easy induction.
\end{eproof}

\smallskip

We are able now to generalize Proposition \ref{prop:xmoment} and get
uniform moment estimates of every order (under assumptions on the
initial condition).

\begin{ecor}\label{cor:momentT}
  If $W$ satisfies (A) or if $W$ and $V$ satisfy (A') (where $m$ is
  defined), then for all $t>0$ and all $k\geq 1$ there exists a
  constant $c(k)$ such that for all $i$
  $$
  \sup_{t\geq 0} \dE |X_t^{i,N}|^{2k} \leq c(k)
  (1+\mu_0(|x|^{2mk})) \, .
  $$
  If $W$ satisfies (A) a similar result holds with $Y_t^{i,N}$ in
  place of $X_t^{i,N}$.
\end{ecor}
\begin{eproof}
  We write the proof under $(A')$. Let $B_k(s)=\dE(|X_s^{i,N}|^{k})$.
  As in the previous propositions we shall use It\^o's formula up to
  the stopping time $T$ and then go to the limit. Using our hypotheses
  we get that for some nonnegative $\lambda'$
\begin{eqnarray*}
B_{2k}(t) & \leq & B_{2k}(0)  + 2k(C+2k-1) \dE \int_0^t |X_s^{i,N}|^{2k-2} ds - 2k \lambda' \, \dE
\int_0^t |X_s^{i,N}|^{2k} ds  \\
& & - {2k\over N} \, \dE\sum_{j=1}^N \int_0^t \left(\nabla
W(X_s^{i,N}-X_s^{j,N})\right)\cdot X_s^{i,N} |X_s^{i,N}|^{2k-2} ds.
\end{eqnarray*}
To bound the last term above, we use (A3) i.e.
$$
|\nabla W(X_s^{i,N}-X_s^{j,N})|\leq M (1+|X_s^{i,N}-X_s^{j,N}|^m) \, ,
$$
and H\"{o}lder inequality in order to obtain the following upper bound
$$
2^{(2k-1)/2k} \, M \, \int_0^t \,
\left[ \dE   |X_s^{i,N}|^{2k}\right]^{\frac{2k-1}{2k}} \,
\left[ \dE \PAR{1+|X_s^{i,N}-X_s^{j,N}|^{2mk}}\right]^{1/2k} \, ds \, .
$$
Now we may use exchangeability and Proposition \ref{prop:momdiff}
to obtain
\begin{eqnarray*}
B_{2k}(t) & \leq  & B_{2k}(0)
+ 2k(C+2k-1) \, \dE \int_0^t |X_s^{i,N}|^{2k-2} ds
-2k\lambda'\dE\int_0^t|X_s^{i,N}|^{2k} ds
\\ & & + c(k)
(1+\mu_0(|x|^{2mk})^{1/2k} \, \int_0^t \, \left[\dE\PAR{
|X_s^{i,N}|^{2k}}\right]^{\frac{2k-1}{2k}} ds \, .
\end{eqnarray*}
As usual, one can get that $B_{2k}$ is differentiable and satisfies
$$
B_{2k}'(t)\leq c(k)
B_{2k-2}(t)-2kB_{2k}(t)+c(k)(1+\mu_0(|x|^{2mk}))B_{2k}(t)^{\frac{2k-1}{k}}.
$$
Since, for every $\varepsilon>0$, it exists $c$ such that
$|x|^{2k-2} \leq c + \varepsilon |x|^{2k}$ and $a^{\frac{2k-1}{2k}}
\leq 1 + a$ we thus obtain
$$
B_{2k}'(t)  \leq   c(k)(1+\mu_0(|x|^{2mk}))-\tilde \lambda B_{2k}(t),
$$
for some $\tilde \lambda>0$ and we can conclude using Gronwall one
more time.
\end{eproof}

\begin{erem}
  If $V$ is identically 0 we may simply remark that
  $$|X_t^{i,N}|^{k} \leq 3^{k-1} \left(|X_0^{i,N}|^{k}+|\sqrt 2
    B_t^i|^k + \left(\frac 1N \, \sum_{j=1}^N \, \int_0^t C
      (1+|X_s^{i,N}-X_s^{j,N}|^{m}) ds\right)^k \right)
  $$
  and directly conclude with the help of Proposition
  \ref{prop:momdiff}, thanks to the convexity of $x \rightarrow
  |x|^k$.
  One could mimic the preceding proof to get the uniform estimate of
  the moments of the projected particle system under assumptions (A).
\end{erem}

\medskip

\subsection{Identification and existence of solutions of the nonlinear SDE}

We prove here how a tightness criterion may ensure that the empirical
law of the particles converges to the solution of the nonlinear SDE
(\ref{eq:nonlin}) proving thus the existence of such a solution.

\begin{elem}\label{lem:tension}
  If $W$ satisfies (A) or $W$ and $V$ satisfy (A'), then for all $s$
  and $t$ smaller than $T$, all $k\geq 2$ and all $i$
  $$
  \dE |X_t^{i,N}-X_s^{i,N}|^{2k} \leq
  C(k,T)(1+\mu_0(|x|^{m(m+2k-1)})) \, |t-s|^{\frac 32} \, .
  $$
\end{elem}
\begin{eproof}
  Let $0\leq s \leq t \leq T$ and $k\geq 1$.
\begin{eqnarray*}
\dE |X_t^{i,N}-X_s^{i,N}|^{2k} & \leq & 2k(2k-1) \dE \int_s^t |X_u^{i,N}-X_s^{i,N}|^{2k-2} du
\\ &  & + \, (2k) \dE \int_s^t |\nabla V(X_u^{i,N})||X_u^{i,N}-X_s^{i,N}|^{2k-1} du \\
& & + {2k\over N} \, \dE\sum_{j=1}^N \int_s^t \left(|\nabla
W(X_u^{i,N}-X_u^{j,N})|\right)|X_u^{i,N}-X_s^{i,N}|^{2k-1} du .
\end{eqnarray*}
Hence using (A3), Corollary \ref{cor:momentT} and H\"{o}lder (with the
best choice of exponents) we get that for $k\geq 1$
$$
\dE |X_t^{i,N}-X_s^{i,N}|^{2k} \leq
C(k,T)(1+\mu_0(|x|^{m(m+2k-1)})) \, |t-s| \, .$$
Plugging this
estimate into the previous inequality and using Cauchy-Schwarz , we
obtain the desired result.
\end{eproof}

It is well known that Lemma \ref{lem:tension} implies that the sequence of the laws of $(s
\rightarrow X_s^{1,N})_N$ defined on $\mathcal{C}([0,T],dR^d)$ is tight. In order to build a
solution to the nonlinear SDE (\ref{eq:nonlin}) we may now follow some standard routine in
mean-field particle systems. We here follow the one in \cite{Mel} Theorem 4.1.4. Thanks to
Proposition 4.2.2 in \cite{Mel} and to Lemma \ref{lem:tension}, the empirical measures $\pi_N$
defined on $\mathcal P(\mathcal P(\mathcal{C}([0,T],\dR^d)))$ by $\pi_N = 1/N \, \sum_{i=1}^N \,
\delta_{X_.^{i,N}}$ is tight too. According to the end of section 4.2 in \cite{Mel}, for any limit
point $\pi_\infty$, any $Q \in \mathcal P(\mathcal{C}([0,T],\dR^d))$ is $\pi_\infty$ a.s. a
solution to the nonlinear martingale problem (up to time $T$). Actually on one hand the proof here
is simpler since we have no jumps, but on the other hand the drifts are unbounded (but with
polynomial growth) so that to justify passage to the limit we have to use the (uniform in $N$)
moment estimates in Corollary \ref{cor:momentT} and Proposition \ref{prop:momdiff}. Details are
straightforward and left to the reader, who will check that a sufficient condition is
$\mu_0(|x|^{m^2}) < +\infty$.
\smallskip

But consider two solutions $(\bar X_.,\bar Z_.)$ of the nonlinear SDE
(\ref{eq:nonlin}), built with the same Brownian motion and the same
initial condition, and introduce an independent copy $(\bar X'_.,\bar
Z'_.)$. Thanks to Proposition \ref{prop:barmoment} in the next
subsection, we know that all processes have a finite second order
moment. Hence, it holds
\begin{eqnarray*}
A(t) = \dE [|\bar X_t - \bar Z_t|^2] & = & - \, \int_0^t \, 2 \dE [(\nabla V(\bar X_s) - \nabla
V(\bar Z_s)).(\bar X_s - \bar Z_s)] ds \\ & - &  \int_0^t \, 2 \dE [(\nabla W(\bar X_s-\bar X'_s)
- \nabla W(\bar Z_s- \bar Z'_s)).(\bar X_s - \bar Z_s)] ds \, .
\end{eqnarray*}
But again since $\nabla W(-a) = - \, \nabla W(a)$ the later integral
can be rewritten
$$
- \int_0^t \, \dE [(\nabla W(\bar X_s-\bar X'_s) -
\nabla W(\bar Z_s- \bar Z'_s)).\left((\bar X_s-\bar X'_s) - (\bar Z_s
  - \bar Z'_s)\right)] ds \, .
$$

Using \eqref{cond:conv} and the local Lipschitz property we see that
$- (\nabla V(x)-\nabla V(y)).(x-y) \leq \beta |x-y|^2$ a similar
result holding for $W$. Using $(a+b)^2 \leq 2a^2 + 2b^2$ and the fact
that $\bar X - \bar Z$ and $\bar X' - \bar Z'$ have the same law we
get
$$
A(t) \leq 3\beta \, \int_0^t A(s) ds \, ,
$$
so that $A(t)=0$.  Hence we have proved (strong) uniqueness for the
nonlinear SDE. As for linear SDE, this notion of uniqueness implies
uniqueness in law.

\smallskip

As a byproduct, we obtain that $\pi_\infty=\delta_Q$ for some $Q$
which is the unique solution of \eqref{eq:nonlin} and that $\pi_N$
goes to $\delta_Q$.  \medskip

What we have obtained is the following: there exists an unique
probability measure $Q$ defined on $C([0,T],\dR^d)$ such that for all
smooth $f$
$$
(\omega,t) \mapsto f(\omega_t) - f(\omega_0) - \int_0^t \,
\left(\Delta f(\omega_s) - \nabla V\cdot\nabla f(\omega_s) - \nabla W*Q_s(\omega_s)\cdot\nabla f(\omega_s) \right) ds
$$
is a $Q$ martingale, with
$$
\nabla W*Q_s(\omega_s)=\int \nabla W(\omega_s - y) Q_s(dy).
$$

Let $H(x)=\int \nabla W(x - y) Q_s(dy)$. Thanks to Corollary
\ref{cor:momentT}, $H(x) \leq C (1+|x|^m)$ (provided
$\mu_0(|y|^{m^2})$ is finite) and is local Lipschitz (thanks to (A3)).
Thus the SDE $$dZ_t =\sqrt 2 dB_t - \nabla V(Z_t) \, dt \, - \, H(Z_t)
\, dt$$
has a strongly unique solution up to its explosion time. Since
$Q$ is a solution, it is the only one. The solution is thus non
explosive, and Girsanov theory tells us that $Q$ is absolutely
continuous w.r.t. the Wiener measure, provided the drift is of finite
energy (see e.g. \cite{CL94} Proposition 2.3), i.e.
$$
\int_0^T \, \int \, \left(|\nabla V|^2 + |H|^2\right)(x) \, Q_s(dx)
\, ds \, < \, +\infty \, .
$$
This condition is satisfied provided $\mu_0\left(|x|^{2m^2}\right) < +\infty$.
Thus for all $s>0$, $Q_s$ is absolutely continuous w.r.t. Lebesgue
measure, with a density $u_s$.  \medskip

Let us summarize our results.
\begin{ethm}\label{thm:exist}
  Assume that $W$ satisfies (A) or that $W$ and $V$ satisfy (A').
  Assume in addition that $\mu_0(|x|^a) < +\infty$ for $a \, = \, \max
  (m(m+3),2m^2)$.

  Then the nonlinear SDE \eqref{eq:nonlin} has an unique strong
  solution $Q$.

  Furthermore for all $t>0$, $Q_t$ (the law of $Q$ at time $t$) is
  absolutely continuous w.r.t.  Lebesgue measure, with density $u_t$
  satisfying for all $T>0$,
  $$
  \sup_{0<t\leq T} \int |x|^{2k} u_t(x) dx \leq C(T)
  (1+\mu_0(|x|^{2mk})) \, ,
  $$
  and $t \rightarrow u_t$ is a solution of \eqref{eq:edp}.

  Finally $t \rightarrow u_t dx$ is the unique solution of
  \eqref{eq:edp} among the set of continuous flows of measures $t \to
  \nu_t$ satisfying for all $T>0$,
  $$
  \int_0^T \, \int \, |x|^{2m^2} d\nu_t dt < +\infty \, .
  $$
\end{ethm}

The only thing it remains to prove is the last statement. Let $t
\rightarrow \nu_t$ be a solution of \eqref{eq:edp} satisfying the
(finite energy) condition above. Then, according to Theorem 4.18 in
\cite{CL94}, a solution of the (linear time inhomogeneous) SDE $$dZ_t
= \sqrt 2 dB_t - \nabla V (Z_t) dt - (\nabla W * \nu_t)(Z_t) dt $$
with initial law $\mu_0$ exists and furthermore its law at time $t$ is
given by $\nu_t$. Hence $Z_.$ is a solution of the nonlinear SDE.
Uniqueness of the later implies that $\nu_t=u_t dx$.  \bigskip

{\bf About the literature.} If existence and uniqueness have been
extensively discussed in the framework of mean field interacting
particle systems with bounded interactions, the case of unbounded
interactions was not much studied. The above proof shows that
convexity at infinity allows us to essentially mimic the bounded case,
just changing the tightness criterion to be used, and using some
recent aspects of stochastic calculus related to singular diffusion
processes.

If $d=1$ another version of Theorem \ref{thm:exist} is obtained in
\cite{brtv98} Theorem 3.1 following a completely different way. One
can notice that some moment condition on $\mu_0$ similar to ours is
also required therein. Our condition (as well as the one in
\cite{brtv98}) is certainly non sharp. If a large part of the method
in \cite{brtv98} can be extended to the $d$ dimensional case, some
aspects (strongly using monotonicity) require additional work.

For the nonlinear PDE, for $d=1$ and $m=2$ an existence and
uniqueness result is stated in \cite{BCCP} p.983. However the initial
measure has to be absolutely continuous with a $\mathcal{C}^2$ density
satisfying $\int |x|^4 d\mu_0 < +\infty$ (while we need a 10 instead
of a 4), and uniqueness holds for classical ($\mathcal{C}^2$) solutions.

The conclusion is that such results are certainly not useless. Since
we shall need stronger moment assumptions in the sequel, we did not
try to obtain the sharpest conditions in Theorem \ref{thm:exist} (the
interested reader may indeed remark that one need for the use of lemma
\ref{lem:tension} for tension only a dependence in time of the order
$|t-s|^{1+\epsilon}$ so that the condition on the initial measure may
be weakened).  \bigskip

\subsection{Uniform moment control for the nonlinear SDE}

We start with the proof of the moment control for any solution of the
nonlinear SDE.

\begin{eprop}\label{prop:barmoment}
  Assume that $W$ satisfies (A) or that $W$ and $V$ satisfy (A'). Then
  $$
  \sup_{t\geq 0} \dE|\bar X_t|^{2} \le \mu_0(|x|^2) + K \, .
  $$
\end{eprop}
\begin{eproof}
  Let $s \rightarrow \bar X'_s$ an independent copy of $s \rightarrow
  \bar X_s$. Then
\begin{eqnarray*}
\dE(|\bar X_t|^{2})&=&\dE(|\bar X_0|^{2})-2k\int_0^t\dE\left( \bar X_s \cdot(\nabla V(\bar X_s)+\nabla
W*u_s(\bar X_s))\right)ds + 2d t\\
&=&\dE(|\bar X_0|^2)-2\int_0^t\dE(\bar X_s\cdot \nabla V(\bar X_s))ds- \int_0^t\dE((\bar X_s-\bar
X'_s)\cdot\nabla W(\bar X_s-\bar X_s'))ds+ 2 dt
\end{eqnarray*}
so that if we denote $v(t)=\dE(|\bar X_t|^2)$, we may first
differentiate w.r.t. time as we did before, and then use the
hypotheses. Remark that, since $W$ is convex at infinity and
symmetric,
$$
x\cdot\nabla W(x)\ge x\cdot\nabla W(0) + \lambda_W |x|^2 - C \geq
-C.
$$
It follows
$$
v'(t)\le -2\lambda_V v(t)+(2d+3C) - 2\lambda_W \left(v(t) - \dE^2(\bar
  X_t)\right)\, .
$$
If (A') holds, $\lambda_V > 0$ and Gronwall's lemma concludes the
proof. If $V=0$, $\lambda_W >0$ and we have assumed that $E(\bar
X_t)=0$ so that Gronwall's lemma concludes the proof.
\end{eproof}

Now, if we assume that existence and uniqueness hold for the non
linear SDE, the results in the previous two subsections imply some
control of the moment of the solution, just taking limits. But the
previous a priori bound was necessary to complete the proof of Theorem
\ref{thm:exist}.

\section{Propagation of chaos}

This section is devoted to the comparison of the behavior, for a fixed
number of particles, of the difference between one particle and the
solution of the non-linear SDE. Ideally, it will be uniform on time
and will decrease quickly to $0$ as the number of particles increases.
We will however see that to get such an estimate, we have to introduce
a new convexity assumption, ensuring strict convexity except at some
points.  \smallskip

All results in this section are written under assumption (A) i.e. with
a vanishing confinement potential $V$. Replacing $Y$ by $X$, the same
results hold under assumption (A'), modifying the proofs in the same
way as we did for various statements in the previous section.
\medskip

\subsection{A first control on the error}

We prove here a first uniform control on the mean square error between
one particle and the solution of the nonlinear SDE. In the sequel we
shall always make assumptions ensuring existence and uniqueness of
strong solutions. Hence we may build solutions for \eqref{eq:projpart}
(or \eqref{eq:part}) and \eqref{eq:nonlin} with the same Brownian
motions and the same initial random variables (obtaining thus and
i.i.d. sample $(\bar X_.^i)_{i=1,...,N}$ of $\bar X_.$).

\begin{ethm}\label{thm:nonunifbound}
  Assume that $W$ satisfies (A). Suppose that the law $\mu_0$ admits a
  large enough polynomial moment. Then there exists $K>0$ such that
\begin{equation}
\label{eq:propchaos}
\sup_{t\ge0}\dE\left(|Y^{i,N}_t-\bar X^i_t|^2\right)\le K.
\end{equation}
\end{ethm}

This result is of course not so good as it doesn't imply the so called
uniform propagation of chaos, since the last bound does not tend to 0
as $N$ tends to infinity. However under the condition that $W$ is the
sum of an uniformly convex function and a Lipschitz compactly
supported one, one may show the following bound, for some positive $a$
$$
\dE\left(|Y^{i,N}_t-\bar X^i_t|^2\right)\le K e^{at}/N
$$
which gives us the non uniform propagation of chaos. The proof of
such an estimate is the conjunction of the proof below and standard
estimates for bounded Lipschitz drifts for mean field particle systems
(see again \cite{Mel}).

\begin{eproof}
  Introduce
  $$
  \bar Y^{i,N}_t=\bar X^i_t-\frac{1}{N}\sum_{j=1}^N\bar X^j_t .
  $$
  Then
\begin{eqnarray*}
\dE\left(|Y^{i,N}_t-\bar X^i_t|^2\right)&\le&2\dE\left(|Y^{i,N}_t-\bar Y^i_t|^2\right)+
2\dE\left(|\bar Y^{i,N}_t-\bar X^i_t|^2\right)\\
&=&2\dE\left(|Y^{i,N}_t-\bar Y^i_t|^2\right)+\frac{2}{N}\dE(|\bar X^i_t|^2)
\end{eqnarray*}
since $\dE(\bar X_t^i)=0$ and the $\bar X_t^j$'s are independent. The
second term is bounded by some $K/N$ according to Proposition
\ref{prop:barmoment} so that we focus on the first term. But
\begin{eqnarray*}
Y^{i,N}_t-\bar Y^{i,N}_t&=&-\frac{1}{N}\int_0^t\sum_{j=1}^N\left[\nabla W(Y^{i,N}_s-Y^{j,N}_s)-
\nabla W*u_s(\bar X^i_s)\right]ds\\
&&-\frac{1}{N}\int_0^t\sum_{j=1}^N\nabla W*u_s(\bar X^j_s)ds
\end{eqnarray*}
so that
$$
\sum_{i=1}^N |Y^{i,N}_t-\bar
Y^{i,N}_t|^2=-\frac{2}{N}\sum_{i,j=1}^N\int_0^t
(A_{ij}(s)+B_{ij}(s)+C_{ij}(s)) ds
$$
with
\begin{eqnarray*}
A_{ij}(s)&=&(\nabla W(Y^{i,N}_s-Y^{j,N}_s)-\nabla W(\bar Y^{i,N}_s-\bar Y^{j,N}_s))\cdot(Y^{i,N}_s-\bar Y^{i,N}_s),\\
B_{ij}(s)&=&(\nabla W(\bar X^{i}_s-\bar X^{j}_s)-\nabla W*u_s(\bar X^i_s))\cdot(Y^{i,N}_s-\bar Y^{i,N}_s),\\
C_{ij}(s)&=&\nabla W*u_s(\bar X^i_s))\cdot(Y^{i,N}_s-\bar Y^{i,N}_s).
\end{eqnarray*}
Let us deal with the first term,
\begin{eqnarray*}
&&\sum_{i,j=1}^N A_{ij}(s) =\frac{1}{2}\sum_{i,j=1}^N (A_{ij}(s)+A_{ji}(s))\\
&&\qquad =\frac{1}{2}\sum_{i,j=1}^N(\nabla W(Y^{i,N}_s-Y^{j,N}_s)-\nabla W(\bar Y^{i,N}_s-\bar Y^{j,N}_s)).
((Y^{i,N}_s-Y^{j,N}_s)-(\bar Y^{i,N}_s-\bar Y^{j,N}_s))\\
&&\qquad \ge \lambda N\sum_{i=1}^N |Y^{i,N}_s-\bar Y^{i,N}_s|^2-CN^2/2
\end{eqnarray*}
using \eqref{cond:conv} and properties of vectors in ${\cal M}$.

Remark now that
\begin{eqnarray*}
-\dE \sum_{j=1}^NB_{ij}(s) &\le& \left(\dE|Y^{i,N}_s-\bar Y^{i,N}_s|^2\right)^{1/2}
\left(\dE|\sum_{j=1}^N(\nabla W(\bar X^i_s-\bar X^j_s)-\nabla W*u_s(\bar X^i_s)|^2\right)^{1/2}\\
&\le&c\sqrt{N}\left(\dE|Y^{i,N}_s-\bar Y^{i,N}_s|^2\right)^{1/2}
\end{eqnarray*}
by Cauchy-Schwarz inequality, the polynomial growth of $\nabla W$, the
controls of moments previously established and the key remark
$$
\dE(\nabla W(\bar X^i_s-\bar X^j_s)-\nabla W*u_s(\bar X^i_s))=0.
$$
The last term involving $C_{ij}(s)$ is controlled using the same
tools. We finally get, defining
$$
\alpha(s)=\dE(|Y^{i,N}_s-\bar Y^{i,N}_s|^2)
$$
and reasoning as we did in the previous section, the differential
inequality
$$
\alpha'(s)\le-2\lambda
\alpha(s)+\frac{c}{\sqrt{N}}\sqrt{\alpha(s)}+C/2 .
$$
Now we may use $\sqrt {\alpha(s)}/\sqrt N \le 1/2(\varepsilon^{-1}
\frac 1N + \varepsilon \alpha(s))$ and choose
$\varepsilon=2\lambda/c$. The previous inequality becomes a classical
$$
\alpha'(s)\le -\lambda \alpha(s)+(c/\varepsilon N) +C/2
$$
for which we can use Gronwall's lemma.
\end{eproof}

\textbf{Remark.} \quad When $C=0$ we get the uniform (in time)
propagation of chaos result obtained in \cite{malrieu03} with rate $1/
N$. Of course in our case (convexity at infinity only) such an uniform
propagation of chaos cannot be expected (more precisely cannot be
shown with usual tools). Indeed, the collision potential $W$ obliges
the particle to mainly stay in the same region, where there is no more
attraction between particles.  \smallskip

\subsection{Uniform propagation of chaos}

In view of Theorem \ref{thm:nonunifbound} and its proof, we have to
reinforce the convexity assumption in order to prove the uniform
propagation of chaos phenomenon.

The condition \eqref{supercond} in the introduction and inspired by the work of
Carrillo-McCann-Villani \cite{CMV2} is recalled below:

we say that condition $\mathbf {C(A,\alpha)}$ holds if there exist $A,\alpha>0$ such that for any
$0<\epsilon<1$,
\begin{equation*}
\forall x,y\in\dR^d,\qquad (x-y). (\nabla W(x)-\nabla W(y))\ge A\epsilon^\alpha(
|x-y|^2-\epsilon^2).
\end{equation*}

We can now prove

\begin{ethm}\label{th:upropchaos}
Assume that $W$ satisfies $\mathbf {C(A,\alpha)}$ and (A). Suppose  that the law $\mu_0$ has a
large enough polynomial moment. Then there exists $K>0$ such that
\begin{equation}
\label{eq:propchaos}
\sup_{t\ge0}\dE\left(|Y^{i,N}_t-\bar X^i_t|^2\right)\le \frac{K}{N^{1\over 1+\alpha}}.
\end{equation}
\end{ethm}

\begin{eproof}
The proof follows the same lines than the one of Theorem \ref{thm:nonunifbound}. First recall that
\begin{eqnarray*}
\dE\left(|Y^{i,N}_t-\bar X^i_t|^2\right)&\le& 2\dE\left(|Y^{i,N}_t-\bar
Y^i_t|^2\right)+\frac{2}{N}\dE(|\bar X^i_t|^2)
\end{eqnarray*}
and the second term is of a better order thanks to Proposition \ref{prop:barmoment} so that we
focus on the first term. For this term the only modification is the control of $A_{ij}(s)$ where
we replace convexity at infinity by \eqref{supercond}. This yields $$ \sum_{i,j=1}^N A_{ij}(s) \ge
A\epsilon^\alpha\left( N\sum_{i=1}^N |Y^{i,N}_s-\bar Y^{i,N}_s|^2-\epsilon^2 N^2/2\right) ,$$ so
that the differential inequality satisfied by $\alpha(t)$ becomes
$$\alpha'(s)\le-2A\epsilon^\alpha(\alpha(s)-\epsilon^2)+\frac{c}{\sqrt{N}}\sqrt{\alpha(s)}.$$
Applying Theorem \ref{thm:nonunifbound} we know that $\alpha(s)\leq K$ for some $K>1$.

We may take $\epsilon=\sqrt{\alpha(s)}/2\sqrt K <1$ and get
 $$\alpha'(s)\le- J \,  \alpha(s)^{1+\alpha/2}+\frac{c}{\sqrt{N}}\sqrt{\alpha(s)}$$ with
 $J=\frac{2A}{(2 \sqrt K)^{\alpha}} \, (1 - \frac{1}{4K})$. Define $\beta(s)=\sqrt {\alpha (s)}$.
 Then $$\beta'(s) + (J/2) \beta^{1+\alpha}(s) \le \frac{c}{2 \sqrt{N}}$$ so that $$\beta(s) \le
 C/N^{1/2(1+\alpha)}$$ for any $s$ such that $\beta'(s) \ge 0$. Since $\beta(0)=0$ it easily
 follows that $\beta(s) \le
 C/N^{1/2(1+\alpha)}$ everywhere, hence the result.
\end{eproof}

\textbf{Remark.} \quad If (A') holds we have to assume that $V$ satisfies $\mathbf {C(A,\alpha)}$
and $W$ is convex.

\section{Ergodicity of the nonlinear SDE}

\subsection{Drift condition and existence of a stationary measure}

In fact, as long as we consider only the particle system, the condition (\ref{cond:conv}) is
sufficient to ensure the ergodicity, and even exponential ergodicity, of the particle system.
Indeed we have that $(X_t^{i,N})_{1\le i\le N}$ is strongly Feller aperiodic and we can use
Down-Meyn-Tweedie's drift condition \cite{DMT}. This condition is the following: there exists some
$\Phi\ge1$ with compact level sets such that for some positive $\lambda,b$
$${\cal L} \Phi \le -\lambda \Phi +b$$
where ${\cal L}$ is the generator of the particle system. It is trivially verified here with the
function $\Phi(x)=\|x\|^2+1$ according to (\ref{cond:conv}). It of course implies the existence of
an invariant measure $\pi^N$ and that there exists $\delta(N)>0$ such that denoting $P^N_t$ the
semigroup associated with the particle system
$$ \| P^N_t(x,\cdot)-\pi^N(\cdot)\|_{TV}\le c \, e^{- t \delta(N)}  \, \Phi(x).$$
As a byproduct, we get the existence of an invariant measure for the nonlinear SDE (hence an
equilibrium for the nonlinear PDE), provided we have uniform (in time) propagation of chaos.
Indeed the previous inequality together with uniform propagation of chaos show that the family
$(u_t)_{t>0}$ is a Cauchy family in $\mathbf L^1(dx)$, which is a complete space.

Such a result is however not useful to control the convergence to equilibrium for the nonlinear
SDE (or even the nonlinear PDE). Indeed there is no close form for $\delta(N)$ which heavily
depends on $N$, so that even the uniform propagation of chaos cannot give good estimations for the
rate of convergence of $(\bar X_t)$ to its invariant measure.  Note that even a strict convexity
condition does not seem to be useful for such an approach.

\subsection{Logarithmic Sobolev inequality}

Under a global uniform convexity condition (i.e. $C=0$ in
(\ref{cond:conv})), Malrieu proves in \cite{malrieu03} that the
invariant measure of the centered particles system i.e.
$$
u^N_\infty(y) dy =Z_N^{-1}\exp\left(-{1\over
    2N}\sum_{i,j=1}^NW(y_i-y_j)\right) \, dy
$$
where $dy$ is Lebesgue measure on the hyperplane $\mathcal M$,
satisfies a logarithmic Sobolev inequality with a constant $1/\lambda$
independent of the dimension using the famous Bakry-\'Emery criterion.
So that if $u^N_t$ stands for the law of $(Y^{i,N}_t)_{1\le i\le N}$
we have that
$$
\text{Ent}(u^N_t|u^N_\infty)\le e^{-2\lambda t}
\text{Ent}(u^N_0|u^N_\infty).
$$
Projecting on the first coordinate, using a $T_2$ inequality and,
on one hand the dimension free exponential convergence, and on the
other hand the uniform (in time) propagation of chaos, he essentially
proves the following

$$
W_2(u_t,u_\infty)\le  C \, e^{-\lambda t}
\sqrt{\text{Ent}(u_0|u_\infty)}.
$$

This result was obtained by Carrillo-McCann-Villani \cite{CMV}, and
improved in Carrillo-McCann-Villani \cite{CMV2} where the authors
replace the square root of the initial entropy by $W_2(u_0,u_\infty)$
(recall that a log-Sobolev inequality implies a $T_2$ inequality, i.e.
$W_2(u_0,u_\infty) \le C \sqrt {\text{Ent}(u_0|u_\infty)}$).

It is of course tempting to use the same approach in the non strictly
convex case.

The first point is that under assumption $\mathbf{C(A,\alpha)}$, we
may no longer use Bakry-\'Emery condition for the measure $u^N_\infty$
as the Hessian may degenerate. We are thus obliged to use some
perturbation argument.

To fix ideas and to simplify the presentation, let us restrict
ourselves to a polynomial potential, namely $W(x)=|x|^{2+\alpha}$.
Consider first a crude attempt: let $\beta(x)=ax^2/(1+x^2)$ with some
positive $a$: $\beta$ is bounded, and has positive curvature $2a$ at
0. Choose now $a=1/N$ and consider the perturbed measure
$$\nu^N_\infty(dy)=Z'_N \exp\left(\frac 1N \,
  \sum_{i,j}\beta(|y_i-y_j|)\right)u^N_\infty(dy).$$
By arguments
developed in Malrieu \cite{malrieu03} one may then easily establish
that $\nu^N_\infty$ has positive curvature bounded below by $1/N$.
Thus it satisfies a logarithmic Sobolev inequality with constant $eN$
remarking that the perturbation potential has oscillations bounded by
one and applying the standard Holley-Stroock perturbation argument. We
may then use the entropic convergence associated to such a logarithmic
Sobolev inequality, the $T_2$ inequality projected on the first
coordinate and uniform propagation of chaos (Theorem
\ref{th:upropchaos}) to get
$$
W_2(u_t,u_\infty)\le \left(\frac{K}{N^{1\over 1+\alpha}}\right)^{1/2}
+e^{- t/4eN} \sqrt{2e N \, \text{Ent}(u_0|u_\infty)}.
$$
Optimizing in $N$ ensures that
$$
W_2(u_t,u_\infty)\le \frac{K'\log^{1\over 2(1+\alpha)} t}{t^{1\over
    2(1+\alpha)}}
$$
which is quite a bad estimate of the speed since
Carrillo-McCann-Villani obtained in such a case a speed
$t^{-(1/\alpha)}$ (we shall recover this speed in the next
subsection).  \smallskip

Note finally that this approach via logarithmic Sobolev inequality has
two major drawbacks: on one hand the right hand side depends on the
entropy of the initial data, on the other hand it cannot be used to
compare the behavior of two solutions starting from different initial
measures.  We will see in the next section how to solve both these
problems by a direct simple approach.

\subsection{A direct control of the $\mathbf{L^2}$-Wasserstein distance}

The goal of this section is to show that even if the logarithmic Sobolev inequality cannot give
useful result, one may use a direct approach to obtain a (uniform in $N$) control of the
Wasserstein distance between two solutions of the particle system starting from different points.
So, combining this with the uniform in time propagation of chaos, we recover in an elegant and
simple way the results of Carrillo-McCann-Villani \cite{CMV2}.

\begin{ethm}\label{th:w2}
Assume that $W$ satisfies (A) and the convexity condition ${\mathbf C}(A,\alpha)$. Let $u_t$ and
$v_t$ be the unique solutions of the nonlinear PDE with initial conditions respectively $u_0$ and
$v_0$. We assume for simplicity that both $u_0$ and $v_0$ have an exponential moment (or a large
enough polynomial moment in order to ensure existence and uniqueness).

Then $t \mapsto W_2^2(u_t,v_t)$ is non-increasing. Furthermore there exists $t_1\leq
(2^{2+\alpha}/3) \, \log (W^2_2(u_0,v_0))/A$ such that
\begin{eqnarray}\label{w2particle}
W_2^2(u_t,v_t)& \le & e^{-(3A/2^{2+\alpha}) t} \, W_2^2(u_0,v_0) \quad \textrm{ if } t<t_1 \, ,
\\ & \le & \left(1+A \, \left(\alpha/(2+\alpha)\right)^{1+\alpha/2}(t-t_1)\right)^{-2/\alpha}
\quad \textrm{ if } t\ge t_1 \, . \nonumber
\end{eqnarray}
In particular, if $t> (1+\eta) (2^{2+\alpha}/3) \, \log (W^2_2(u_0,v_0))/A$ for some $\eta>0$, one
has
$$W_2^2(u_t,v_t) \le \left(1+(A/(1+\eta)) \,
\left(\alpha/(2+\alpha)\right)^{1+\alpha/2} \, t\right)^{-2/\alpha}.$$ In addition for all $t\ge
0$ one has $$W_2^2(u_t,v_t) \le \left(W_2^{-\alpha}(u_0,v_0)+A \,
\left(\alpha/(2+\alpha)\right)^{1+\alpha/2} \, t\right)^{-2/\alpha}.$$
\smallskip

Under an uniform convexity condition for $W$ (i.e. ${\mathbf C}(A,0)$), the convergence is
exponential i.e.
$$
W_2(u_t,v_t)\le C e^{-At}W_2(u_0,v_0).
$$
\end{ethm}
Note that the final assertion completes the result of Malrieu
\cite{malrieu03} who cannot compare via logarithmic Sobolev
inequalities the distance between two solutions.

Remark also that the bound obtained for the control in Wasserstein
distance in the convergence towards the stationary measure is also
better as it is now $W_2(u_0,u_\infty)$ which controls the decay and
which is smaller (since a $T_2$ inequality holds) than $C
\sqrt{\text{Ent}(u_0|u_\infty)}$.

This will be explained later.

As will be seen from the proof, there is no use here of heavy (but
sharp) technology as optimal transport in length space as developed by
Carrillo-McCann-Villani \cite{CMV2} and thus we believe that such an
approach will give good results for other models (reinforced
diffusion,...).

\begin{eproof}
  In what follows we shall always consider coupling consisting in
  picking the same Brownian motions $(B^i)$ for the two particle
  systems, and only choose an ad-hoc coupling for the initial random
  variables.

  Denote $u^{1,N}_t$ the law of the first particle of the centered
  particle system ($Y$) starting with initial law $u_0$ and consider
  $v^{1,N}_t$ the law of the first particle of the centered particle
  system starting with initial law $v_0$. Let us first remark that for
  each coupling $g_0$ of $u_0$ and $v_0$
  $$
  W_2^2(u^{1,N}_t,v^{1,N}_t)\le {1\over N} \,
  \dE_{g_0}\left(\sum_{i=1}^N|Y^{i,N}_t-Y'^{i,N}_t|^2\right)
  $$
  where $(Y^{i,N})$ is given by the centered particle system where
  each particle starts with law $u_0$ , $(Y'^{i,N})$ with measure
  $v_0$ and the initial law of the pairs $(Y^{i,N}_0,Y'^{i,N}_0)$ are
  given by independent copies of $g_0$. The subscript for the
  expectation is related to the initial law.

  We then get by It\^o's formula and symmetry (starting now from
  points $y$ and $y'$)
\begin{eqnarray*}
&&\dE_{y,y'}\left(\sum_{i=1}^N|Y^{i,N}_t-Y'^{i,N}_t|^2\right)\\
&&\quad =\sum_{i=1}^N|y_i-y'_i|^2-{1\over N}\sum_{i,j} \dE_{y,y'} \int_0^t\langle
Y^{i,N}_s-Y'^{i,N}_s~,~\nabla W(Y^{i,N}_s-Y^{j,N}_s)-\nabla W(Y'^{i,N}_s-Y'^{j,N}_s)
\rangle ds\\
&&\quad=\sum_{i=1}^N|y_i-y'_i|^2-{1\over 2N}\sum_{i,j}\dE_{y,y'} \int_0^t \xi_{ij}(s)ds,
\end{eqnarray*}
with
$$
\xi_{ij}(s)=
\langle
(Y^{i,N}_s-Y^{j,N}_s)-(Y'^{i,N}_s-Y'^{j,N}_s)~,~\nabla W(Y^{i,N}_s-Y^{j,N}_s)-\nabla
W(Y'^{i,N}_s-Y'^{j,N}_s)\rangle.
$$
We may now differentiate in time and then use condition
$\mathbf{C}(A,\alpha)$ to get
\begin{eqnarray*}
{d\over dt}\dE_{y,y'}\left(\sum_{i=1}^N|Y^{i,N}_t-Y'^{i,N}_t|^2\right)
&=&-{A\epsilon^\alpha\over 2N}\sum_{i,j}\dE_{y,y'}(|(Y^{i,N}_t-Y^{j,N}_t)-(Y'^{i,N}_t-Y'^{j,N}_t)|^2-\epsilon^2)\\
&=&-A\epsilon^\alpha \, \dE_{y,y'}\sum_i(|Y^{i,N}_t-Y'^{i,N}_t|^2-\epsilon^2)
\end{eqnarray*}
which gives us by exchangeability
$$
{d\over dt}\dE_{y,y'}\PAR{|Y^{i,N}_t-Y'^{i,N}_t|^2}
\le-A\epsilon^\alpha \, \dE(|Y^{i,N}_t-Y'^{i,N}_t|^2-\epsilon^2).
$$
Denote $\xi(t)=\dE_{y,y'}\PAR{|Y^{i,N}_t-Y'^{i,N}_t|^2}$, we then obtain
for all $\epsilon<1$
$$
\xi'(t)\le -A\epsilon^\alpha(\xi(t)-\epsilon^2).
$$
This inequality for $\varepsilon$ going to $0$, implies $\xi'(t)\le
0$, i.e. $\xi$ is non-increasing. As a byproduct we get that $t
\mapsto W_2(u_t^{1,N},v_t^{1,N})$ is also non-increasing. Indeed if we
choose $g_0$ as the optimal coupling for the quadratic cost, it holds
$$
W_2^2(u_0,v_0)=\xi(0) \ge \xi(t) \ge W^2_2(u_t^{1,N},v_t^{1,N})\, ,
$$
and we get the result just shifting the initial time.  \smallskip

Now we can separate two cases: either $\xi(t)> 1$ or $\xi(t) \le 1$.
Note that there exists some $t_1\ge 0$ such that the first case holds
for $t<t_1$ and the second one for $t\ge t_1$.

If $t<t_1$ we may choose any $\varepsilon$, for instance here we
choose $\varepsilon=1/2$ and obtain (since $\xi(t)>1$, $\varepsilon^2
=1/4 < \xi(t)/4$)
$$
\xi'(t)\le -A(\alpha) \xi(t)
$$
with $A(\alpha)=(3A/4) \, (1/2)^{\alpha}$, which gives by
Gronwall's lemma
$$
\xi(t)\le e^{-A(\alpha) t} \, \xi(0).
$$
Choosing again the optimal coupling $g_0$ we obtain the first part
of the result and the fact that $t_1 \le \log (W^2_2(u_0,v_0))/A(\alpha)$.

For the second part, if $t\ge t_1$ choose $\epsilon^2=\alpha
\xi(t)/(\alpha+2)$ to get
$$
\xi'(t)\le -A\left({\alpha\over 2+\alpha}\right)^{\alpha/2}{2\over
  2+\alpha}~\xi(t)^{1+\alpha/2}.
$$
Integrating this differential inequality we get
$$
\xi(t)\le \left(1+B(\alpha)(t-t_1)\right)^{-2/\alpha},
$$
with $B(\alpha)=A \, \left(\alpha/(2+\alpha)\right)^{1+\alpha/2}$.
Our choice of $\varepsilon$ gives the optimal constant (for this
method).  \smallskip

Of course writing $\mathbf{C}(A,\alpha)$ as we did is a little bit
artificial, and if we want a more homogeneous estimate we may remark
that $\mathbf{C}(A,\alpha)$ implies
$$
(x-y)\cdot(\nabla W(x)-\nabla W(y))\ge
A(\epsilon/W_2(u_0,v_0))^\alpha( |x-y|^2-\epsilon^2)
$$
for $\epsilon < W_2(u_0,v_0)$. So we may always choose
$\epsilon^2=\alpha \xi(t)/(\alpha+2)$ and get for all $t$,
$$
\xi(t)\le
\left(W_2^{-\alpha}(u_0,v_0)+B(\alpha)t\right)^{-2/\alpha}.
$$

We use then the uniform propagation of chaos property to transfer the
inequality from the particles to the solutions of the nonlinear SDE:
$$
W_2(u_t,v_t)\le
W_2(u_t,u^{1,N}_t)+W_2(u^{1,N}_t,v^{1,N}_t)+W_2(v_t,v^{1,N}_t),
$$
and take the previous (uniform in $N$) estimation for the middle
term and the uniform in time estimation for the first and third term
and let $N$ go to infinity. Of course if $t_1$ depends on $N$, its
bound does not, so that we may find a converging subsequence and get
the result.
\end{eproof}

\textbf{Remarks.} \quad We have seen in subsection 4.1 that there
exists an invariant measure $u_\infty$ for the nonlinear PDE. The
bound for $W_2(u_t,v_t)$ obtained by introducing the particle system
of order $N$, and then choosing choosing $v_0=u_\infty^N$ allows us to
prove this existence too (using completeness of the Wasserstein
distance). In addition the decay of the Wasserstein distance trivially
implies uniqueness of $u_\infty$, at least in the set of measures
having some large enough polynomial moment.  \smallskip

Hence, we recover in the previous result the asymptotic rate of
convergence to equilibrium obtained by Carrillo-McCann-Villani. But it
seems that the result improves upon theirs for small times where we
are able to describe some exponential decay. However the value of
$t_1$ is not explicit, so that for practical issues, this initial
exponential decay is not really tractable.  \smallskip

We choose here to develop the point of view of the particles; prove
the decay of the Wasserstein distance for the particles system and
then transfer this decay to the solution of the nonlinear SDE via
uniform propagation of chaos. It is however easy to develop the same
line of proof directly for the control of the Wasserstein distance of
two solutions of the nonlinear SDE.  \smallskip

Note finally that we did not consider here coefficient diffusion other
than constant, as granular media equation is formulated with a
constant one. It would however not be difficult to introduce a
condition enabling to obtain the same decay. Indeed consider
$$dZ_t=b(Z_t)dt+\sigma(Z_t)dW_t$$
where $W_t$ is the usual Wiener
process in $\dR^d$, and assume that for all $0<\epsilon<1$ there
exists positive $A$ and $\alpha$ such that
$$
(x-y)\cdot(b(x)-b(y))+{1\over 2}tr(
(\sigma(x)-\sigma(y))(\sigma(x)-\sigma(y))^t)\le
-A\epsilon^\alpha(\|x-y\|^2-\epsilon^2)
$$
then we easily derive the same polynomial decay for the Wasserstein
distance of two solutions $Z_t$ and $\tilde Z_t$ with different
initial conditions. Note that this type of condition hold with
possible degeneracy of the diffusion coefficient.  \medskip

Again in this section we should assume that (A') holds, $V$ satisfies
$\mathbf {C(A,\alpha)}$ and $W$ is convex.


\section{Concentration inequality}

The main goal of this section is to complete the results on
convergence of the particle system and of the nonlinear system by
providing a (non asymptotic) deviation inequality. This inequality
will be written for additive functionals of the particles and then
allows us to estimate integrals with respect to the stationary measure
of the nonlinear PDE. Once gain it is interesting to get uniform in
times estimation to be able to simulate at fixed time the particles
and use them for the evaluation with no loss at each time of the
constant in the concentration. Therefore we first prove in the general
framework an uniform $T_1$ inequality for solution of SDE subject to
some convexity at infinity condition. We then show how to use them on
our example.

\subsection{Uniform transportation cost inequality under convex at infinity condition}

Let $X$ be the solution of the following stochastic differential
equation:
\begin{equation}
  \label{eq:eds}
dX_t = dB_t+b(X_t)dt,
\end{equation}
where $B$ is a standard Brownian motion on $\dR^d$ and $b$ a smooth
function from $\dR^d$ to $\dR^d$.

In the case when $b(x)=x$, $X$ is the well-known Ornstein-Uhlenbeck
process and
$$
\sup_{t\geq 0}\dE\SBRA{e^{\de\ABS{X^x_t-Y^x_t}^2}} < +\infty,
$$
(where $X^x$ and $Y^x$ are two independent copies starting both at
$x$) is finite if and and only if $\de <1/2$, since the law of $X_t$
is the Gaussian measure
$$
\mathcal Law (X_t^x)=\mathcal N\left(xe^{-t},\frac{1-e^{-2t}}{2}\right).
$$
Our aim is to extend this assertions to the case when the drift is
confining only outside a compact set, and for non-constant diffusion
coefficient, so that $X$ is the solution of
\begin{equation}
  \label{eq:eds2}
dX_t =\sigma(t) dB_t+b(X_t)dt.
\end{equation}
We use here the formalism of transportation cost inequalities in $W_1$
distance for which a practical criterion based on the integrability of
the exponential of the square of the distance is sufficient \cite{DGW}
(and Bolley-Villani \cite{BV} or Gozlan \cite{gozlan} for a better
evaluation of the constant) and which implies interesting deviation
inequalities.
\begin{eprop}\label{t1}
  Suppose that there exist $\la,A >0$ and $C$ such that, for every
  $x,y\in\dR^d$,
\begin{equation}\label{eq:reconv}
(x-y)\cdot(b(x)-b(y))\leq -\la \ABS{x-y}^2+ C, \qquad \|\sigma\|_{HS}\le A
\end{equation}
then for any $\de<\la/2A$ and any $x\in\dR^d$,
$$
\sup_{t\geq 0}\dE\SBRA{e^{\de\ABS{X^x_t-Y^x_t}^2}} \leq 1 + \PAR{Ad+C+1}e^{\de (Ad+C+1)/(\la -2\de
A)},
$$
where $X^x$ and $Y^x$ are two independent copies of \eqref{eq:eds}
starting at $x$.

So the law of $X_t^x$ satisfies a $T_1$ inequality with a constant
$\mathfrak C$ independent of time and initial position.
\end{eprop}

Recall that a $T_1$ inequality with constant $\mathfrak C$ for a measure $\mu$ reads as:
$$W_1(\nu,\mu) \leq \sqrt {\mathfrak C \, \text{Ent}(\nu|\mu)}$$
for any
$\nu$.

\begin{eproof}
  The process $Z=X-Y$ is a solution of
  $$
  dZ_t=\sigma(X_t)dW^1_t-\sigma(Y_t)dW^2_t+(b(X_t)-b(Y_t))dt,
  $$
  with initial condition $Z_0=0$, where $W^1$ and $W^2$ are two
  independent Brownian motions. Itô's formula ensures that
\begin{eqnarray*}
  e^{\de\ABS{X_t-Y_t}^2}-e^{\de\ABS{X_s-Y_s}^2}&\le&
2\de\int_s^t\!(X_u-Y_u)\cdot(b(X_u)-b(Y_u))e^{\de\ABS{Z_u}^2}du\\
&&+2\de A\int_s^t\!(d+2\de \ABS{X_u-Y_u}^2)e^{\de\ABS{Z_u}^2}du+ M_t-M_s\\
&\leq& 2\de\int_s^t\!(Ad+C+(2A\delta-\la)\ABS{X_u-Y_u}^2)e^{\de\ABS{Z_u}^2}du+ M_t-M_s,
\end{eqnarray*}
where $M$ is a local martingale with quadratic variation given by
$$
\DP{M}_t =
8A\delta^2\int_0^t\!\ABS{X_u-Y_u}^2e^{\de\ABS{X_u-Y_u}^2}\,du.
$$
For every $\de<\la/2A$, let us denote by $K$ the quantity
$$
K=\sqrt{\frac{Ad+C+1}{\la-2\de A}}.
$$
For every $x\in\dR^d$,
\begin{eqnarray*}
(Ad+C-(\la-2\de A)\ABS{x}^2)e^{\de\ABS{x}^2}&\leq&
(Ad+C-(\la-2\de A)\ABS{K}^2)e^{\de\ABS{x}^2}\ind_\PAR{\ABS{x}\geq K}
+(Ad+C)e^{\de K^2}\\
&\leq&-e^{\de \ABS{x}^2}+e^{\de \ABS{x}^2}\ind_\PAR{\ABS{x}\leq K}
+(Ad+C)e^{\de K^2}\\
&\leq&-e^{\de \ABS{x}^2}+(Ad+C+1)e^{\de K^2}.
\end{eqnarray*}
For $R>0$, introduce $T_R$ the first time when $Z$ exits the ball of
radius $R$ and define $\al_R$ the function defined by
$$
\al_R(t)=\dE\SBRA{e^{\de\ABS{Z_{t\wedge T_R}}^2}}.
$$
Then
$$
\al_R(t)-1 \leq 2\de (Ad+C+1)e^{\de K^2}t ,
$$
so that we may let $R$ go to infinity, and have shown that the
exponential moment is finite.

Defining now $\al(t)=\dE\SBRA{e^{\de\ABS{Z_{t}}^2}}$ , one gets, for
$0\leq s\leq t$,
$$
\al(t)-\al(s)\leq 2\de (Ad+C+1)e^{\de K^2}(t-s)-2\de
\int_s^t\!\al(u)\,du.
$$
As a consequence, $\al$ satisfies the following differential
inequality
$$
\al'(t)\leq -2\de\al(t) + 2\de (Ad+C+1)e^{\de K^2}.
$$
Denoting by $\be(t)=\al(t)\exp(2\de t)$, this implies that,
$$
\be'(t)\leq 2\de (Ad+C+1)e^{\de K^2} e^{2\de t},
$$
and then,
$$
\be(t)\leq \be(0)+ (Ad+C+1)e^{\de K^2} \PAR{e^{2\de t}-1}.
$$
Besides, $X$ and $Y$ have the same initial condition so
$\be(0)=\al(0)=1$. As a conclusion, we get for all $t\geq 0$,
$$
\al(t)\leq e^{-2\delta t} + (Ad+C+1)e^{\de K^2},
$$
which achieves the first part of the proof. The final statement
follows from \cite{DGW} Theorem 3.1. An explicit expression of
$\mathfrak C$ is derived in \cite{gozlan} chapter VII or in
Bolley-Villani \cite{BV}.
\end{eproof}

An important consequence of a transportation inequality is that we
easily obtain deviation inequality for Lipschitz functions. Indeed for
all Lipschitz functions $F$ with $\|F\|_{Lip}\le 1$ and all positive
$r$
$$\dP\left(F(X^x_t)-\dE(F(X^x_t))\ge r\right)\le  e^{-r^2/{\mathfrak C}}.$$
\medskip

Remark that the previous proof extends to the case when the initial
law $\mu_0$ satisfies $$\int \, e^{\delta |x-y|^2} \, \mu_0(dx) \,
\mu_0(dy) \, < \, + \infty$$
(that is $\mu_0$ satisfies a $T_1$
inequality) just choosing two independent variables $X_0$ and $Y_0$ of
law $\mu_0$.

\subsection{A concentration inequality for the stationary measure}

We may now recombine results inherited from the previous sections to
get an useful inequality for the evaluation of $\int f du_\infty$ when
$f$ is a Lipschitz function. From the previous section, if (A) holds
and $\mu_0$ satisfies a $T_1$ inequality, the particle system
satisfies a $T_1$ inequality with constant ${\mathfrak C}N$ (for some
${\mathfrak C}$ independent of time), which thus leads to the
following:

for $f$ Lipschitz (in $\dR^d$) with $\|f\|_{Lip}\le 1$, for all
positive $r$ and all $t$
$$
\dP\left({1\over N}\sum_{k=1}^N f(X^{k,N}_t)-\dE f(X^{1,N}_t)\ge
  r\right)\le e^{-Nr^2/{\mathfrak C}}.
$$
Remark now that for a Lipschitz function satisfying $\|f\|_{Lip}\le
1$, if the convexity condition $\mathbf {C(A,\alpha)}$ holds, the
uniform propagation of chaos of Theorem \ref{th:upropchaos} tells us
$$
|\dE f(X^{1,N}_t)-\int f(y) u_t(y) dy|\le W_1(u^{1,N}_t,u_t) \le
W_2(u^{1,N}_t,u_t)\le \left(K\over N^{1\over 1+\alpha}\right)^{1/2}
$$
so that for all $r \ge \left(K\over N^{1\over
    1+\alpha}\right)^{1/2}$
$$
\dP\left({1\over N}\sum_{k=1}^N f(X^{k,N}_t)-\int f(y) u_t(y) dy\ge
  r-\left(K\over N^{1\over 1+\alpha}\right)^{1/2}\right) \le
e^{-Nr^2/{\mathfrak C}}.
$$
We may now use convergence in $W_2$ distance of the solutions of
the nonlinear SDE towards the stationary measure given in Theorem
\ref{th:w2} to get

\begin{eprop} For all $r\ge \left(K\over N^{1\over
1+\alpha}\right)^{1/2}+\sqrt{\beta(t)}$,
\begin{equation}
\dP\left({1\over N}\sum_{k=1}^N f(X^{k,N}_t)-\int f(y) u_\infty(y) dy\ge r-\left(K\over N^{1\over
1+\alpha}\right)^{1/2}- \sqrt{\beta(t)}\right)\le e^{-Nr^2/{\mathfrak C}},
\end{equation}
where $\beta(t)$ is one of the functions governing the decay of
$W^2_2(u_t,u_\infty)$ described in Theorem \ref{th:w2}.
\end{eprop}
\medskip

It is quite hard to imagine to extend to our case uniform result over
Lipschitz function (i.e.  deviation of the $W_1$ distance between the
empirical law of the particles and the stationary measure) as in
Bolley-Guillin-Villani \cite{BGV} as it requires a dynamic coupling
which can be achieved only for potential $W$ whose Hessian is bounded
(in the sense of matrix).

Let us finally note than one can use an Euler-Maryama scheme
preserving square exponential integrability and with good stability
property \cite[Th1 and Th 4]{LMS} to simulate the particle system
leading to the same concentration inequality. Note that the recurrence
property needed for the stability of this adaptive scheme in
\cite{LMS} is exactly our condition of convexity at infinity.



\bigskip
\noindent
P. Cattiaux: \'Ecole Polytechnique, CMAP, CNRS 756, 91128 Palaiseau Cedex FRANCE and
Universit\'e Paris X Nanterre, Equipe MODAL'X, UFR SEGMI, 200 avenue de la
R\'epublique, 92001 Nanterre cedex, FRANCE.

\noindent
Email: cattiaux@cmapx.polytechnique.fr

\medskip\noindent
A. Guillin: CEREMADE, UMR CNRS 7534, Place du Mar\'echal De Lattre De Tassigny
75775 PARIS CEDEX 16 - FRANCE.

\noindent
Email: guillin@ceremade.dauphine.fr\\
Web:  http://www.ceremade.dauphine.fr/\raisebox{-4pt}{$\!\widetilde{\phantom{x}}$}guillin/\\

\medskip\noindent
F. Malrieu: IRMAR, Universit\'e Rennes 1, Campus de Baulieu, 35042 Rennes cedex, France.

\noindent
Email: florent.malrieu@univ-rennes1.fr\\
Web:  http://name.math.univ-rennes1.fr/florent.malrieu/\\

\end{document}

%% file: macros.tex
\newcommand{\dE}{\ensuremath{\mathbb{E}}}

\newcommand{\dN}{\ensuremath{\mathbb{N}}}

\newcommand{\dP}{\ensuremath{\mathbb{P}}}

\newcommand{\dR}{\ensuremath{\mathbb{R}}}

\newtheorem{ethm}{Theorem}[section]

\newtheorem{ecor}[ethm]{Corollary}

\newtheorem{eprop}[ethm]{Proposition}

\newtheorem{elem}[ethm]{Lemma}

\newtheorem{edefi}[ethm]{Definition}

\newtheorem{erem}[ethm]{Remark}

\newcommand{\proofend}{~$\rhd$}
\newcommand{\proofbegin}{~$\lhd$}

\newenvironment{eproof}
               {\noindent {\emph{\textbf{Proof}}}\\\proofbegin~}
               {\proofend\\}

\newcommand{\ABS}[1]{\ensuremath{{\left| #1 \right|}}} 
\newcommand{\PAR}[1]{\ensuremath{{\left(#1\right)}}} 
\newcommand{\SBRA}[1]{\ensuremath{{\left[#1\right]}}} 
\newcommand{\DP}[1]{\ensuremath{{\left<#1\right>}}} 
\newcommand{\PD}[2]{\ensuremath{{\frac{\partial #1}{\partial #2}}}} 

\renewcommand{\phi}{\varphi}

\renewcommand{\geq}{\geqslant}

\newcommand{\ind}{\mathrm{1}\hskip -3.2pt \mathrm{I}}

\newcommand{\al}{\alpha}

\newcommand{\be}{\beta}

\newcommand{\de}{\delta}

\newcommand{\la}{\lambda}

